\documentclass[12pt,a4paper]{article} 
\usepackage{amsmath,pstricks,latexsym,theorem,epsfig} 
\usepackage{amsfonts,amssymb,latexsym,epsfig,amsmath}
\usepackage{amsmath,latexsym,theorem} 
\newtheorem{Theorem}{Theorem}[section] 
\newtheorem{Definition}{Definition}[section] 
\newtheorem{Proposition}{Proposition}[section] 
\newtheorem{Lemma}{Lemma}[section] 
\newtheorem{Corollary}{Corollary}[section] 
 
\theorembodyfont{\rmfamily}

\newcommand{\rec}[1]{{(\ref{#1})}} 
\newcommand{\be}{\begin{equation}}
\newcommand{\ee}{\end{equation}}
\newcommand{\R}{\mathbb{R}}
\newcommand{\N}{\mathbb{N}}

\newcommand{\Z}{\mathbb{Z}}

\def\ti{\tilde}
\def\lf{\left}
\def\rg{\right}

\def\ds{\displaystyle}

\def\Om{\Omega}

\def \N{I\!\!N} 

\def \R{I\!\!R}

\def \Z{Z\!\!\!Z} 
 
\def\11{1\!\!1}

 \def\res{\mathop{\hbox{\vrule height 7pt width .5pt 
depth 0pt\vrule height .5pt width 6pt depth 0pt}}\nolimits}
\def\ds{\displaystyle}

\def\rec#1{{(\ref{#1})}}

\newcommand{\ba}{\begin{array}} 
\newcommand{\ea}{\end{array}}

\begin{document} 
 \date{}
 \title{\sc Fractional Harmonic Maps into Manifolds in odd dimension $n>1$}
\author{ Francesca Da Lio\thanks{Department of Mathematics, ETH Z\"urich, R\"amistrasse 101, 8092 Z\"urich, Switzerland.} \thanks{Dipartimento di Matematica Pura ed Applicata, Universit\`a degli Studi di Padova. Via Trieste 63, 35121,Padova, Italy.}     }
\maketitle 
\begin{abstract}
 In this paper we consider critical points of the following nonlocal energy
\begin{equation} 
{\cal{L}}_n(u)=\int_{\R^n}| ( {-\Delta})^{n/4} u(x)|^2 dx\,,
\end{equation}
where $u\colon H^{n/2}(\R^n)\to{\cal{N}}\,$ ${\cal{N}}\subset\R^m$ is a compact $k$ dimensional smooth manifold without boundary and   $n>1$ is an odd integer.
Such critical points are called $n/2$-harmonic maps into ${\cal{N}}$.
We prove that $\Delta ^{n/2} u\in L^p_{loc}(\R^n)$ for every $p\ge 1$ and thus 
$u\in C^{0,\alpha}_{loc}(\R^n)\,.$ The local H\"older continuity of $n/2$-harmonic maps   is based on regularity results obtained in \cite{DL1} for nonlocal
Schr\"odinger systems with an antisymmetric potential   and on suitable  {\it 3-terms commutators } estimates.
\end{abstract}
 {\small {\bf Key words.} Harmonic maps, nonlinear elliptic PDE's, regularity of solutions, commutator estimates.}\par
 {\small { \bf  MSC 2000.}  58E20, 35J20, 35B65, 35J60, 35S99}

\section{Introduction}
In the paper \cite{DLR2} the authors considered $1/2$-harmonic maps in $\R$ with values in a k-dimensional sub-manifold ${\cal{N}}\subset\R^m$ ($m\ge 1$), which is  smooth, compact and without boundary.
We recall that $1/2$-harmonic maps are functions $u$ in the space
  $
 \dot{H}^{1/2}(\R,{\cal{N}})=\{u\in \dot{H}^{1/2}(\R,\R^m):~~u(x)\in {\cal{N}}, {\rm a.e},\}\,,
 $
which are critical points for perturbation of the type $\Pi^{N}_ {\cal{N}}(u+t\varphi)$, ($\varphi\in C^{\infty}$ and
$\Pi^{N}_ {\cal{N}}$ is the normal projection on $ {\cal{N}}$)
 of the functional
\begin{equation}\label{lagr1}
{\cal{L}}_1(u)=\int_{\R}| ( {-\Delta})^{1/4} u(x)|^2 dx\,,
\end{equation}
  (see Definition 1.1 in \cite{DLR1})\,.
  The operator $ ( {-\Delta})^{1/4} $ on ${\R}$ is defined by means of the the Fourier transform as follows  $$\widehat{ ( {-\Delta})^{1/4} u}=|\xi|^{1/2}\hat u\,,$$
 (given a function f, $\hat f$ or ${\cal{F}}$ denotes the Fourier transform of $f$).\par

The Lagrangian \rec{lagr1}  is invariant with respect to the M\"obius group and it satisfies 
  the following identity
 \begin{eqnarray*}
\int_{\R}| ( {-\Delta})^{1/4} u(x)|^2 dx=\inf\left\{\int_{\R^2_+}|\nabla \tilde u|^2 dx:~\tilde u\in W^{1,2}(\R^2,\R^m),~~\mbox{trace $\tilde u=u$}\right\}\,.
\end{eqnarray*}

  The Euler Lagrange equation associated to the nonlinear problem \rec{lagr1} can be written as follows :
 
 \be
 \label{zz6}
  ( {-\Delta})^{1/2}u\wedge \nu(u)=0\quad\quad\mbox{ in }\cal{D}'({\R})\quad,\
 \ee
where $\nu(z)$ is the Gauss Map at $z\in\cal{N}$ taking values into the grassmannian $\ti{G}r_{m-k}({\R}^m)$ of oriented $m-k$ planes in ${\R}^m$, which to every
point  $z\in\cal{N}$ assigns the unit $m-k$ vector defining   
  the oriented normal $m-k-$plane to $T_z\cal{N}$\,.
 The $C^{0,\alpha}_{loc}$ regularity of $1/2$ harmonic maps was deduced from a  key  result
 obtained in \cite{DLR2} concerning with nonlocal
 linear Schr\"odinger system in $\R$ with an antisymmetric potential of the type:
  \begin{equation}
 \label{zz1}
\forall i=1\cdots m\quad\quad\quad ( {-\Delta})^{1/4} v^i =\sum_{j=1}^m\Omega_j^i\, v^j\,,
\end{equation}
where $v=(v_1,\cdots,v_m)\in L^2({\R},{\R}^m)$ and $\Omega=(\Omega_i^j)_{i,j=1\cdots m}\in L^2({\R},so(m))$ is an $L^2$ maps from ${\R}$ into the space $so(m)$ of $m\times m$ antisymmetric matrices.

 It is   natural to extend the above mentioned results   to $n/2$ harmonic maps in $\R^n$,  with values in a k-dimensional sub-manifold ${\cal{N}}\subset\R^m$, where $m\ge 1 $ and $n=2p+1$ is an odd integer.
 By analogy with the case $n=1$, $n/2$ harmonic maps  are functions $u$ in the space
  $
 \dot{H}^{n/2}(\R^n,{\cal{N}})=\{u\in \dot{H}^{n/2}(\R^n,\R^m):~~u(x)\in {\cal{N}}, {\rm a.e},\}\,,
 $
which are critical points for perturbation of the type $\Pi^{N}_ {\cal{N}}(u+t\varphi)$, ($\varphi\in C^{\infty}(\R^n,\R^m))$  of the functional
 \begin{equation}\label{lagrangian}
{\cal{L}}_n(u)=\int_{\R^n}| ( {-\Delta})^{n/4} u(x)|^2 dx\,.
\end{equation}
The Euler Lagrange equation associated to the non linear problem \rec{lagrangian} can be written as follows :
 
 \be
 \label{zz7}
  ( {-\Delta})^{n/2}u\wedge \nu(u)=0\quad\quad\mbox{ in }\cal{D}'({\R})\,.
 \ee
 The  Euler Lagrange in the form (\ref{zz7}) is hiding fundamental properties such as it's elliptic nature...etc and is difficult to use directly 
for solving problems related to regularity and compactness. One of the first task is then to rewrite it in a form that will make some of it's
analysis features more apparent. This is the purpose of the next proposition. Before to state it we need some additional notations\,.
 
 \medskip
 
Denote by $P^T(z)$ and $P^N(z)$  the projections respectively to the tangent space $T_z\cal{N}$ and to the normal space $N_z\cal{N}$ to  ${\cal{N}}$ at $z\in \cal{N}$ . For $u\in \dot{H}^{n/2}({\R^n},\cal{N})$ we denote simply by $P^T$ and $P^N$ the compositions $P^T\circ u$ and $P^N\circ u$. Under  the assumption that $\cal{N}$ is  smooth, $P^T\circ u$ and $P^N\circ u$ are matrix valued maps in $\dot{H}^{n/2}({\R^n},M_m({\R}))$. 
We will prove the following    useful formulation of the $n/2$-harmonic map equation.

 \begin{Proposition}\label{EulEq}
Let $u\in\dot{H}^{n/2}(\R^n,{\cal{N}} )$ be a weak $n/2$-harmonic map. Then the following equation holds
\be\label{Euler}
 ( {-\Delta})^{n/4}v=\Om\, v+\ti{\Om}_1\, v+\ti{\Om}_2
\ee
where $v\in L^2({\R^n},{\R}^{2m})$ is given by
\[
v:=\left(\begin{array}{l}
P^T ( {-\Delta})^{n/4} u\\
 P^N ( {-\Delta})^{n/4} u\end{array}\right)
\]
and where ${\cal{R}}$ is the Fourier multiplier of symbol $\sigma(\xi)=i\frac{\xi}{|{\xi}|}$. 
\noindent $\Om\in L^2({\R^n},so(2m))$ is given by
$$
\Omega=2\left(\begin{array}{cc}
- \omega &  \omega \\[5mm]
 \omega   &  -\omega  \end{array}\right)\,,$$ the map $\omega$  is in $L^2({\R^n},so(m))$ and given   by
$$
\omega=\frac{ ( {-\Delta})^{n/4} P^T P^T- P^T ( {-\Delta})^{n/4} P^T}{2}\quad.
$$
 Finally   $\ti{\Om}_1=\ti{\Om}_1(P^T,P^N,( {-\Delta})^{n/4} u)$ is  in $L^{(2,1)}({\R^n},M_{2m}({\R}))$ and  $\ti{\Om}_2=\ti{\Om}_2(P^T,P^N )$ is in $\dot W^{-n/2,(2,\infty)}({\R^n},{\R}^{2m})$,    and satisfy
\begin{equation}\label{tildeomega2intr}
\|\tilde\Omega_1\|_{ L^{(2,1)}(\R^n)}\le C( \|P^N\|^2_{\dot H^{n/2}(\R^n)}+ \|P^T\|^2_{\dot H^{n/2}(\R^n)})\,;
\end{equation}
   \begin{equation}\label{tildeomega1intr}
\|\tilde\Omega_2\|_{ \dot W^{-n/2,(2,\infty)}(\R^n)}\le C( \|P^N\|_{\dot H^{n/2}(\R^n)}+ \|P^T\|_{\dot H^{n/2}(\R^n)})\| ( {-\Delta})^{n/4} u\|_{L^{(2,\infty)}(\R^n)}\,.
\end{equation}
\end{Proposition}\par
\medskip
The explicit formulations of $\ti\Omega_1$ and $\ti{\Om}_2$  are given in  Section \ref{harmonic}.   
The control  on $\ti\Omega_1$ and $\ti{\Om}_2$ is a consequence of regularity by compensation results on some operators that we now introduce.
\par

Given $Q\in {\cal{S}}^{\prime}(\R^n, {\cal{M}}_{\ell\times m}(\R^n))$ $\ell\ge 0$\footnote{ ${\cal{M}}_{\ell\times m}(\R)$    denotes, as usual, the space of $\ell\times m$ real matrices.} and $u\in {\cal{S}}^{\prime}(\R^n,\R^m)$, let us define  the   operator $T_n$ as follows.
\begin{equation}\label{opT}
T_n(Q,u)=  ( {-\Delta})^{n/4}[( ( {-\Delta})^{n/4} Q)\,   u]-Q ( {-\Delta})^{n/2} u+ ( {-\Delta})^{n/4} [  Q\,  ( ( {-\Delta})^{n/4} u)]\,.
\end{equation}
  We prove the following commutator estimate. 
  \begin{Theorem}\label{comm1}
Let $u\in\dot  W^{n/2,(2,\infty)}(\R^n)$, $ Q\in \dot H^{n/2}(\R^n)$. Then 
$
T_n(Q,u) \in  \dot H^{-n/2}(\R^n)\,$
 and
\begin{equation}\label{commest1}
||T_n(Q,u) ||_{ \dot H^{-n/2}(\R^n)}\le C||Q||_{\dot{H}^{n/2}(\R^n)}|| ( {-\Delta})^{n/4} u||_{L^{(2,\infty)}(\R^n)}\,.~~~\Box\end{equation}
\end{Theorem}
Theorem \ref{comm1} is a straightforward consequence of the following estimate   for the dual operator of $T_n$ defined by
\begin{equation}\label{optildeT}
 T_n^*(Q,u)=  ( {-\Delta})^{n/4}[( ( {-\Delta})^{n/4} Q)\,   u]- ( {-\Delta})^{n/2}[Q\, u]+ ( {-\Delta})^{n/4} [  Q\,  ( ( {-\Delta})^{n/4} u)]\,.
\end{equation}  
\begin{Theorem}\label{comm2}
Let $u,Q\in \dot H^{n/2}(\R^n)$. Then 
$
 T_n^*(Q,u) \in W^{-n/2,(2,1)}(\R^n)\,,$
 and
\begin{equation}\label{commest2}
|| T_n^*(Q,u)||_{\dot  W^{-n/2,(2,1)}(\R^n)}\le C\|Q\|_{\dot{H}^{n/2}(\R^n)}\|u\|_{\dot{H}^{n/2}(\R^n)}\,.~~~\Box\end{equation}
\end{Theorem}\par

We recall that the spaces $ W^{n/2,(2,\infty)}(\R^n)$ and $W^{-n/2,(2,1)}(\R^n)$ are defined as
$$\dot W^{n/2,(2,\infty)}(\R^n):=\{f\in{\cal{S}}^{\prime}:~~|\xi|^{n}{\cal{F}}[v]\in L^{(2,\infty)}(\R^n)\}\,;$$
$$\dot  W^{-n/2,(2,1)}(\R^n):=\{f\in{\cal{S}}^{\prime}:~~|\xi|^{-n}{\cal{F}}[v]\in L^{(2,1)}(\R^n)\}\,.$$\par
Moreover $\dot W^{n/2,(2,\infty)}(\R^n)$ is the dual of $\dot  W^{-n/2,(2,1)}(\R^n)$.
We refer the reader to   Section \ref{defnot} for the definition of Lorentz spaces $L^{(p,q)}$, $1\leq p,q\le+\infty$ and of the fractional Sobolev spaces\,.\par

Theorems \ref{comm1} and \ref{comm2} correspond respectively to    Theorem 1.2  and Theorem 1.4 in \cite{DLR1} for $n=1$.
The difference with the case $n=1$ is that here we are able to show that $ T_n^*$ is not necessarily in the Hardy space  ${\cal{H}}^1(\R^n)$ but in the bigger space $ W^{-n/2,(2,1)}(\R^n)$. \par
The main result of this paper is the following
\begin{Theorem}\label{regn2hm}
Let  ${\cal{N}}$ be a smooth closed submanifold of ${\R}^m$. Let  $u\in \dot{H}^{n/2}(\R^n,{\cal{N}})$ be a weak $n/2-$harmonic map into ${\cal{N}}$, then  $u\in C^{0,\alpha}_{loc}(\R^n,{\cal{N}})\,.$ \hfill$\Box$
\end{Theorem}
Finally a classical ''elliptic type'' bootstrap argument leads to the following result  
(see \cite{DL1} for the details of this argument).
\begin{Theorem}
\label{th-I.3}
Let ${\cal{N}}$ be a smooth closed submanifold of ${\R}^m$. Let $u$ be a weak $n/2$-harmonic map in $\dot{H}^{n/2}({\R^n},{\cal{N}}))$, then $u$ is $C^\infty$ .
\hfill $\Box$
\end{Theorem}
We mention that Theorem \ref{regn2hm} follows by a slight perturbation of the following Theorem
concerning with the sub-criticality of non-local Schr\"odinger systems in dimension $n>1$. 
The proof of this result is given in    \cite{DL1} and   extends Theorem 1.1 in \cite{DLR2}.
   \begin{Theorem}\label{regschrn}
Let $\Omega\in L^2(\R^n,so(m))$ and $v\in L^2(\R^n)$ be a weak solution of
\begin{equation}\label{eqschbis}
 ( {-\Delta})^{n/4} v=\Omega \,v\,.
\end{equation}
  Then $v\in L^p_{loc}(\R^n)$ for every $1\le p<+\infty$.
\end{Theorem}
We conclude with some comments.\par
  The proof of Theorem \ref{comm1} is not a mere extension of Theorem 1.3 in \cite{DLR2}.
  The fact that we are dealing with the dimension $n>1$, it requires a different  analysis when we split the operator $T_n$ in the so-called paraproducts (see Appendix  \ref{commutators}).\par
  Moreover we mention that in the case of $n=1$ the pseudo-differential operators
  $\nabla$ and $ ( {-\Delta})^{1/2}$ are of the same orders and this permits to write the
  equations for $P^T ( {-\Delta})^{1/4} u$ and $P^N ( {-\Delta})^{1/4} u$ in a similar way (see Section 5 in \cite{DLR2}).
  More precisely $P^T ( {-\Delta})^{1/4} u$ and $P^N ( {-\Delta})^{1/4} u$  satisfy
  \begin{equation}\label{euler1intr}
 ( {-\Delta})^{1/4} ( P^T  ( {-\Delta})^{1/4}u)=T_1(P^T,u)-( ( {-\Delta})^{1/4} P^T) ( {-\Delta})^{1/4} u\,,
\end{equation}
and
\begin{equation}\label{eqstructintr}
  ( {-\Delta})^{1/4} {\cal{R}}[P^N  ( {-\Delta})^{1/4}u)]=S_1(P^N,u)- [ ( {-\Delta})^{1/4} P^N] {\cal{R}} [({-\Delta})^{1/4} u]\,,
\end{equation}
where for $Q,u\in {\cal{S}}^{\prime}(\R^n)$
\[
 S_1(Q,u):= ( {-\Delta})^{1/4}[Q ( {-\Delta})^{1/4} u]-{\cal{R}}  [Q\nabla u]+ [(-\Delta)^{1/4} Q]{\cal{R}}[(-\Delta )^{1/4} u]\,.
 \]
 To write down \rec{euler1intr} and \rec{eqstructintr} we use respectively the fact that $P^T ( {-\Delta})^{1/2}u=0$ and $P^N\nabla u=0\,.$
  In the case   $n>1$,  $\nabla$ and $ ( {-\Delta})^{n/2}$ are   pseudo-differential operators of order respectively $1$ and $n$. The equation for $P^T ( {-\Delta})^{n/4} u$ is similar to equation $(58)$ in
   \cite{DLR2} ($T_1$ is replaced by $T_n$). On the contrary 
   we cannot replace 
    the equation \rec{eqstructintr} by an equation of the form
  \begin{equation}\label{eqstructintrbis}
  ( {-\Delta})^{n/4} {\cal{R}}[P^N  ( {-\Delta})^{n/4}u)]=S_n(P^N,u)- [ ( {-\Delta})^{n/4} P^N] {\cal{R}} ({-\Delta})^{n/4} u]\,,
\end{equation}
where for $Q,u\in {\cal{S}}^{\prime}(\R^n)$
\[
 S_n(Q,u):= ( {-\Delta})^{1/4}[Q ( {-\Delta})^{n/4} u]-{\cal{R}} ( {-\Delta})^{\frac{n-1}{2} }  [Q\nabla u]+ [(-\Delta) ^{n/4} Q] {\cal{R}}[(-\Delta)^{n/4} u]\,.
 \]
 Actually even if  $S_n$ seems the  natural extension of $S_1$,  it does not satisfy the same  regularity estimates as $S_1$, (see \cite{DLR1})\,.\par
 In the case of $n>1$,  the structure equation  
  becomes 
  \begin{equation}\label{eqstructintr3}
 ( {-\Delta})^{n/4} (P^N  ( {-\Delta})^{n/4}u)= ( ( {-\Delta})^{n/4} \bar  {\cal{R}}) f(P^N,u) \,.
\end{equation}
where 
\begin{equation}\label{fintr}
f(P^N,u):=    {\cal{R}}(P^N  ( {-\Delta})^{n/4}u)-  ( {-\Delta})^{ \frac{n}{4}-\frac{1}{2}}  [P^N  \nabla u] \,.
\end{equation}
 We show that the right hand side of equation \rec{eqstructintr3} is in  $W^{-n/2,(2,\infty)}(\R^n)$ and
 $$
 \| (( {-\Delta})^{n/4}\bar  {\cal{R}}) f(P^N,u)\|_{W^{-n/2,(2,\infty)}(\R^n)}\lesssim\|P^N\|_{H^{n/2}}\| ( {-\Delta})^{n/4} u\|_{L^{(2,\infty)}}\,.
 $$
 We conclude by recalling existing results in the literature on regularity of critical points of nonlocal Lagrangians and we refer the reader to \cite{Riv} and \cite{Riv3}  for a complete overview of analogous results in the local case.
 
 The regularity of $1/2$-harmonic maps with values into a sphere has ben first investigated in \cite{DLR1} where new ``three terms comutators" estimates have been obtained by using the technique of paraproducts.
 Analogous results have been extended in \cite{DLR2} to   $1/2$-harmonic maps with values into   general submanifolds.
 
In \cite{Moser} the author considers critical points to the functional
that assigns to any $u\in \dot{H}^{1/2}({\R},{\cal{N}})$ the minimal Dirichlet energy among all possible extensions \underbar{ in $\cal{N}$},
while in the papers \cite{DLR1,DLR2}  the classical $\dot{H}^{1/2}$ Lagrangian corresponds to the minimal Dirichlet energy among all possible extensions \underbar{ in ${\R}^m$}.
Hence the approach in \cite{Moser} consists in working with an \underbar{intrinsic} version of $H^{1/2}-$energy instead of  an \underbar{extrinsic} one.
The drawback of considering the \underbar{intrinsic} energy is that the Euler Lagrange equation is almost impossible to write explicitly and is then \underbar{implicit} .  However the intrinsic version of the $1/2-$harmonic map is more closely related to the existing regularity theory of Dirichlet Energy minimizing maps into $\cal{N}$.
  Finally the regularity of $n/2$ harmonic maps   in odd dimension $n>1$ with values into a sphere has been recently investigated     by Schikorra \cite{Schi} .   In this paper the author extends the results
  obtained in \cite{DLR1} by using an approach based  on compensation arguments introduced by Tartar \cite{Tart85}. 
  
  \medskip

The paper is organized as follows.

  - In Section \ref{defnot} we recall some basic definitions and notations.\par
  - In Section \ref{harmonic} we derive the Euler- Lagrangian equation \rec{Euler} associated to the energy \rec{lagrangian} and we prove Theorem \ref{regn2hm}.\par
  - In Appendix \ref{commutators} we prove the commutator estimates that are used
  in Section \ref{harmonic}.
  
  \section{Preliminaries: function spaces and the fractional Laplacian}\label{defnot}

In this Section we introduce some  notations and definitions that are used in the paper.
\par

For $n\ge 1$,  we denote respectively by ${\cal{S}}(\R^n)$ and  ${\cal{S}}^{\prime}(\R^n)$ the spaces of Schwartz functions and tempered distributions. Moreover given a function $v$ we will denote either by
  $\hat v$ or by  ${\mathcal{F}}[v]$ the Fourier Transform of  $v$ :
  $$\hat v(\xi)={\mathcal{F}}[v](\xi)=\int_{\R^n}v(x)e^{-i \langle \xi, x\rangle }\,dx\,.$$
  Throughout the paper we use the convention that $x,y$ denote variables in the space and
  $\xi$, $\zeta$ the variables in the phase\,.
  
  We recall the definition of fractional Sobolev space (see for instance \cite{T3}).\par
  \begin{Definition}\label{fracsob} 
  For a real $s\ge 0$, 
  $$H^{s}(\R^n)=\lf\{v\in L^2(\R^n):~~|\xi|^s{\cal{F}}[v]\in L^2(\R^n)\rg\}\quad$$
  For a real $s<0$,
   $$H^{s}(\R^n)=\lf\{v\in {\cal{S}}^{\prime} (\R^n):~~(1+|\xi|^2)^s{\cal{F}}[v]\in L^2(\R^n)\rg\}\,.$$
   \hfill $\Box$
   \end{Definition}
   It is known that $H^{-s}(\R^n)$ is the dual of $H^{s}(\R^n)\,.$\par
     
For a submanifold  ${\cal{N}}$ of $\R^m$  we can define   
 $$H^{s}({\R}^n,{\cal{N}})=\{u\in H^{s}({\R}^n,\R^m):~~u(x)\in {\cal{N}}, {\rm a.e.}\}\,.
 $$
 Given $q>1$ we also set
 $$
 W^{s,q}(\R^n):=\{v\in L^q({\R}^n):~~|\xi|^s{\cal{F}}[v]\in L^q(\R^n)\}\,.
 $$
 
 \medskip

 We shall make use of the Littlewood-Paley  dyadic  decomposition of unity that we recall here. Such a decomposition can be obtained as follows.
 Let $\phi(\xi)$ be a radial Schwartz function supported in $\{\xi\in{\R}^n:~|\xi|\le 2\}$, which is
 equal to $1$ in $\{\xi\in{\R}^n: ~|\xi|\le 1\}$\,.
 Let $\psi(\xi)$ be the function given by
 $$
 \psi(\xi):=\phi(\xi)-\phi(2\xi)\,.
 $$
 $\psi$ is then a ``bump function'' supported in the annulus $\{\xi\in{\R}^n:~1/2\le |\xi|\le 2\}\,.$
 
 \medskip
 
 Let $\psi_0=\phi$, $\psi_j(\xi)=\psi(2^{-j}\xi)$ for $j\ne 0 \,.$ The functions $\psi_j$, for $j\in\Z$, are supported in  $\{\xi\in{\R}^n:~2^{j-1}\le |\xi|\le 2^{j+1}\}\,.$
 and  realize a dyadic decomposition of the unity :
  $$
  \sum_{j\in\Z}\psi_j(x)=1\,.
  $$
 We denote further $$\phi_j(\xi):=\sum_{k=-\infty}^j\psi_k(\xi)\,.
 $$ 
 The function $\phi_j$ is supported on  $\{\xi, ~|\xi|\le 2^{j+1}\}$.\par

 We recall the definition of the homogeneous Besov spaces $\dot{B}_{p,q}^s(\R^n)$  and  homogeneous Triebel-Lizorkin spaces $\dot{F}_{pq}^s(\R^n)$ in terms of the above dyadic decomposition.
  \begin{Definition}
  Let $s\in\R$,  $0< p,q\le\infty\,.$ For $f\in{\cal{S}}^\prime (\R^n)$ we set 
   \begin{equation}
   \left.\begin{array}{ll}
 \ds \|u\|_{\dot{B}_{p,q}^s(\R^n)}=\left(\sum_{j=-\infty}^{\infty}2^{jsq}\|{\cal{F}}^{-1}[\psi_j {\cal{F}}[u]]\|_{L^{p}(\R^n)}^q\right)^{1/q}&~~\mbox{if $q<\infty$}\\[5mm]
  \ds  \|u\|_{\dot{B}_{p,q}^s(\R^n)}=\sup_{j\in \Z} 2^{js}\|{\cal{F}}^{-1}[\psi_j {\cal{F}}[u]]\|_{L^{p}(\R^n)}&~~\mbox{if $q=\infty$}\end{array}\right.
    \end{equation}
    When $p,q<\infty$ we also set
   $$ 
   \|u\|_{\dot{F}_{p,q}^s(\R^n)}=\lf\|\left(\sum_{j=-\infty}^{\infty}2^{jsq}|{\cal{F}}^{-1}[\psi_j {\cal{F}}[u]]|^q\right)^{1/q}\rg\|_{L^p}\,.
   $$
   \hfill $\Box$
    \end{Definition}
    The space of all tempered distributions $u$ for which the quantity $\|u\|_{\dot{B}_{p,q}^s(\R^n)}$ is finite is called the homogeneous Besov space with indices 
    $s,p,q$ and it is denoted by $\dot{B}_{p,q}^s(\R^n)$. The space of all tempered distributions $f$ for which the quantity $\|f\|_{\dot F_{p,q}^s(\R^n)}$ is finite is called the homogeneous
    Triebel-Lizorkin space with indices 
    $s,p,q$ and it is denoted by $\dot F_{p,q}^s({\R}^n)\,.$ 
A classical result says \footnote{See for instance \cite{Gra1}} that $\dot{W}^{s,q}(\R^n)=\dot{B}^s_{q,2}(\R^n)=\dot{F}_{q,2}^s(\R^n)$\,.  
    
  Finally we denote    ${\cal{H}}^1(\R^n)$  the homogeneous Hardy Space in $\R^n$. A less classical results \footnote{ See for instance \cite{Gra2}.} asserts that  ${\cal{H}}^1(\R^n)\simeq \dot F^0_{2,1}$ thus we have  
  $$
   \|u\|_{{\cal{H}}^{1}(\R^n)}\simeq \int_{\R^n} \left(\sum_j |{\cal{F}}^{-1}[\psi_j {\cal{F}}[u]]|^2\right)^{1/2}dx\,.
 $$
  \par
  
  We recall that  
 \begin{equation} 
 \dot{H}^{n/2}(\R^n)\hookrightarrow BMO({\R^n})\hookrightarrow \dot{B}^0_{\infty,\infty}(\R^n)\,,\end{equation}
 where $BMO({\R})$ is the space of bounded mean oscillation dual to ${\cal{H}}^1(\R^n)$   (see for instance  \cite{RS}, page 31).
 \par 
 The $s$-fractional  Laplacian of a function  $u\colon\R^n\to\R$ is defined as a pseudo differential operator of symbol $|\xi|^{2s}$ :
 \begin{equation}\label{fract}
 \widehat { ( {-\Delta})^{s}u}(\xi)=|\xi|^{2s} \hat u(\xi)\,.
 \end{equation}
 Finally we introduce the definition of Lorentz spaces (see for instance \cite{Gra1} for a complete presentation of such spaces).
For  $1\le p<+\infty,   1\le q\le +\infty$, the Lorentz space  $L^{(p,q)}(\R^n)$  is the set 
     of measurable functions satisfying
$$
 \left\{\begin{array}{ll}
\int_{0}^{+\infty}(t^{1/p}f^*(t))^{q}\frac{dt}{t}<+\infty, & ~\mbox{if $q<\infty,~p<+\infty$}\\[5mm]
\sup_{t>0}t^{1/p}f^*(t)<\infty & ~\mbox{if $q=\infty,~p<\infty$}\,,\end{array}\right.
 $$
where $f^*$ is the decreasing rearrangement of $|f|\,.$\par
We observe that $L^{p,\infty}(\R^n)$ corresponds to the weak $L^p$ space.
 Moreover for $1<p<+\infty, 1\le q\le\infty$, the space $L^{(\frac{p}{p-1},\frac{q}{q-1})}$ is the dual space of $L^{(p,q)}\,.$
\par
Let us define
$$\dot W^{s,(p,q)}(\R^n)=\{f\in{\cal{S}}^{\prime}:~~|\xi|^s{\cal{F}}[v]\in L^{(p,q)}(\R^n)\}\,.$$
 
 In the sequel we will   often use
  the   H\"older Inequality in the  Lorentz spaces: if $f\in L^{p_1,q_1}, g\in L^{p_2,q_2}$,  with $1\le p_1,p_2,q_1,q_2\le +\infty$. Then
$fg\in L^{r,s},$ with $r^{-1}={p_1}^{-1}+{p_2}^{-1}$ and $s^{-1}={q_1}^{-1}+{q_2}^{-1}\,,$ (see for instance \cite{Gra1}).\par

To conclude we introduce some basic notations.\par
    $B_r(\bar x)$ is the ball of radius $r$ and centered at $\bar x$. If $\bar x=0$ we simply write
  $B_r$\,. If $x,y\in\R^n,$ $x\cdot y$ is the scalar product between $x,y$\,.

 Given a {\em multindex}  
$\alpha=(\alpha_1,\ldots,\alpha_n)$,  where $\alpha_i$ is a
nonegative integer, we denote by  $|\alpha|=\alpha_1+\ldots+\alpha_n\,$ the order of $\alpha$.  \par
   
  For every function $u\colon\R^n\to\R$,  $M(u)$ is the maximal function of $u$, namely
 \begin{equation}\label{maxf}
 M(u)=\sup_{r>0,\ x\in\R^n}|B(x,r)|^{-1}\int_{B(x,r)}|u(y)|dy\,.
 \end{equation}

\section{Euler Equation for $n/2$-Harmonic Maps into Manifolds}\label{harmonic}
We consider a compact $k$ dimensional smooth manifold without boundary ${\cal{N}}\subset\R^m$.   Let  $\Pi_{{\cal{N}}}$ be the 
orthogonal projection on ${\cal{N}}\,.$   Let  $\Pi_{{\cal{N}}}$ be the 
orthogonal projection on ${\cal{N}}\,.$  
 We also consider the   {Dirichlet  energy}  \rec{lagrangian}\,.\par
 The weak ${n/2}$-harmonic maps are defined as critical points
of the functional \rec{lagrangian} with respect to perturbation of the form $\Pi_{{\cal{N}}}(u+t\phi)$, where $\phi$ is
an arbitrary compacted supported smooth map    from $\R $ into $\R^m\,.$
\begin{Definition}\label{weakhalfharm}
We say that $u\in H^{n/2}(\R^n,{\cal{N}})$ is a weak  ${n/2}$-harmonic map if and only if, for every
 maps $\phi\in  H^{n/2}(\R^n,\R^{m})\cap L^{\infty}(\R^n,\R^m)$     we have
\begin{equation}\label{critic}
\frac{d}{dt}{\cal{L}}_n(\Pi_{{\cal{N}}}(u+t\phi))_{|_{t=0}}=0\,.
\end{equation}
\end{Definition}
We introduce some notations.
 We denote by $\bigwedge(\R^m)$ the exterior algebra (or Grassmann Algebra) of $\R^m$ and by the symbol $\wedge$ the {\em exterior or wedge product}. 
For every $p=1,\ldots,m$,   $\bigwedge_p(\R^m)$ is the vector space of $p$-vectors  \par If $(\epsilon_i)_{i=1,\ldots,m}$ is the 
canonical orthonormal  basis of $\R^m$, then every element $v\in \bigwedge_p(\R^m)$ is written as
$v=\sum_{I}v_{I}\epsilon_{I}$ where $I=\{i_1,\ldots,i_p\}$ with $1\le i_1\le\ldots\le i_p\le m$ , $v_I:=v_{i_1,\ldots,i_p} $ and $ \epsilon_{I}=:=\epsilon_{i_1}\wedge\ldots\wedge \epsilon_{i_p}\,.$ \par
By the symbol $\res$ we denote the interior multiplication  $\res\colon \bigwedge_p(\R^m)\times \bigwedge_q(\R^m)\to\bigwedge_{q-p}(\R^m)$ defined as follows. \par
 Let   $\epsilon_I=\epsilon_{i_1}\wedge\ldots\wedge \epsilon_{i_p}$, $\epsilon_J=\epsilon_{j_1}\wedge\ldots\wedge \epsilon_{j_q},$ with $q\ge p\,.$ Then  $\epsilon_I\res \epsilon_J=0$ if $I\not\subset J$, otherwise
  $\epsilon_I\res \epsilon_J=(-1)^M\epsilon_{K}$ where $\epsilon_{K} $ is a $q-p$ vector and $M$ is the number of pairs $(i,j)\in I\times J$ with $j>i\,.$ \par 
 Finally by the symbol $\ast$ we denote the Hodge-star operator, $\ast\colon\bigwedge_p(\R^m)\to\bigwedge_{m-p}(\R^m)$,
defined by $\ast\beta=\beta\res(\epsilon_{1}\wedge\ldots\wedge \epsilon_n) $.  
For an introduction of the Grassmann Algebra we refer the reader to the first Chapter of the book by Federer\cite{fed}\,.\par
 In the sequel we denote by $P^T$ and $P^N$ respectively the tangent and the normal projection to the manifold ${\cal{N}}$.\par
They verify the following properties: $(P^T)^t=P^T,(P^N)^t=P^N$ (namely they are symmetric operators), $(P^T)^T=P^T$,   $(P^N)^N=P^N$,
$P^T+P^N=Id$, $P^NP^T=P^TP^N=0\,.$\par
We set $e=\epsilon_1\wedge\ldots\wedge\epsilon_k$ and $n=\epsilon_{k+1}\wedge\ldots\wedge\epsilon_m$\,. For every $z\in {\cal{N}}$,
$e(z)$ and $n(z)$ give the orientation respectively of the tangent $k$-plane and the normal $m-k$-plane to $T_z{\cal{N}}\,.$\par
We observe that 
   for every $v\in\R^m$ we have
\begin{eqnarray}\label{tangpr}
P^T v&=&(-1)^{m-1}\ast ( ( v\res e) \wedge n) \,.
\end{eqnarray}
\begin{eqnarray}\label{normpr}
P^N v&=&(-1)^{k-1}\ast(e\wedge( v\res n)) \,.
\end{eqnarray}
 Hence $P^N$and $P^T$ can be seen as matrices in $\dot{H}^{n/2}(\R^n,\R^m)\cap
L^{\infty}(\R^n,\R^m)$\,.\par

Next we    write the Euler equation associated to the functional \rec{lagrangian}\,.\par
\begin{Proposition} 
\label{pr-I.1}
All weak $n/2$-harmonic maps $u\in H^{n/2}(\R^n,{\cal{N}})$ satisfy in a weak sense
\par
i)  the equation
\begin{equation}\label{perp}
\int_{\R^n}({\Delta}^{n/2} u )\cdot v \,dx=0,
\end{equation}
for every $v\in \dot{H}^{n/2}(\R^n,\R^m)\cap L^{\infty}(\R^n,\R^m)$ and $v\in T_{u(x)}{\cal{N}}$ almost everywhere, or in a equivalent way\par
ii) the equation 
\begin{equation}\label{wedge}
 P^T{\Delta}^{n/2}u =0~~\mbox{in ${\cal{D}}^\prime\,,$}
 \end{equation}
 or 
\par
iii) the equation
\begin{equation}\label{euler1}
 ( {-\Delta})^{n/4} ( P^T  ( {-\Delta})^{n/4}u)=T_n(P^T,u)-( ( {-\Delta})^{n/4} P^T) ( {-\Delta})^{n/4} u\,,
\end{equation}
where $T_n$ is the operator defined in \rec{opT}\,.
    \end{Proposition}
    The Euler Lagrange equation (\ref{euler1}) can be  completed by the following ''structure equation'':

\begin{Proposition} 
\label{pr-I.2}
All maps in $\dot{H}^{n/2}(\R^n ,{\cal{N}})$ satisfy the following identity
\begin{equation}\label{eqstructbis}
  ( {-\Delta})^{n/4} (P^N  ( {-\Delta})^{n/4}u)= ( ( {-\Delta})^{n/4}\bar  {\cal{R}}) f(P^N,u) \,,
\end{equation}
where 
\begin{equation}\label{f}
 ( ( {-\Delta})^{n/4}\bar  {\cal{R}}) f(P^N,u):=   ( ( {-\Delta})^{n/4}\bar  {\cal{R}}) [  {\cal{R}}(P^N  ( {-\Delta})^{n/4}u)-  ( {-\Delta})^{ \frac{n}{4}-\frac{1}{2}}  [P^N  \nabla u]] 
\end{equation}
 is in $\dot W^{-n/2,(2,\infty)}(\R^n)$ and
\begin{equation}\label{fest}
\| ( ( {-\Delta})^{n/4}\bar  {\cal{R}}) f(P^N,u)\|_{ \dot W^{-n/2,(2,\infty)}(\R^n)}\lesssim \|P^N\|_{\dot H^{n/2}(\R^n)}\| ( {-\Delta})^{n/4}u\|_{L^{(2,\infty)}}\,.
\end{equation}
\end{Proposition}
We give the proof of Proposition \ref{pr-I.2}. For 
  the proof of Proposition \ref{pr-I.1}  we refer the reader to \cite{DLR1}\,.\par\medskip
  
  {\bf Proof of Proposition \ref{pr-I.2}.}
  We first observe that $P^N \nabla u=0$ (see Proposition 1.2 in \cite{DLR1}).
  Thus we can write:
  \begin{eqnarray}\label{f2}
   &&  ( {-\Delta})^{n/4} [P^N  ( {-\Delta})^{n/4}u]= (  ( {-\Delta})^{n/4}\bar {\cal{R}}) {\cal{R}}[P^N  ( {-\Delta})^{n/4}u]\nonumber\\[5mm]
 &&~~~=\underbrace{(  ( {-\Delta})^{n/4} {\cal{R}})(P^N  ( {-\Delta})^{n/4}u)-  ( {-\Delta})^{n/4}  (P^N ( ( {-\Delta})^{n/4}{\cal{R}}u))}_{(1)}\\[5mm]
  &&~~~~+
  \underbrace{  (  ( {-\Delta})^{n/4}\bar {\cal{R}})[P^N ( ( {-\Delta})^{n/4}{\cal{R}}u)]-   (  ( {-\Delta})^{n/4}\bar {\cal{R}})[ ( {-\Delta})^{ \frac{n}{4}-\frac{1}{2}}(P^N \nabla u)]}_{(2)}\nonumber\\[5mm]
  &&~~~= (  ( {-\Delta})^{n/4}\bar {\cal{R}}) f(P^N,u)\,.\nonumber
  \end{eqnarray}
  Corollary  \ref{crwbis} and Theorem \ref{l2} imply respectively that $(1)$ and $(2)\in \dot W^{-n/2,(2,\infty)}(\R^n)$ and 
   $$
  \|(1)\|_{\dot W^{-n/2,(2,\infty)}(\R^n)}\lesssim \|P^N\|_{\dot H^{n/2}(\R^n)}\,\|\Delta^{n/4} u\|_{L^{(2,\infty)}(\R^n)}\,;$$
  $$
  \|(2)\|_{\dot W^{-n/2,(2,\infty)}(\R^n)},\lesssim \|P^N\|_{\dot H^{n/2}(\R^n)}\,\|\Delta^{n/4} u\|_{L^{(2,\infty)}(\R^n)}\,.$$
  
   Hence $( ( {-\Delta})^{n/4}\bar  {\cal{R}}) f(P^N,u)\in  \dot W^{-n/2,(2,\infty)}(\R^n)$ and \rec{fest} holds.~\hfill$\Box$
  \par
  \bigskip
Next we see that by combining  \rec{euler1} and \rec{eqstructbis} we can obtain the  new equation
 \rec{Euler} for the vector field
$v=(P^T  ( {-\Delta})^{n/4}u, P^N  ( {-\Delta})^{n/4}u))$ where an antisymmetric potential appears.  \par
  We introduce the following matrices
\begin{eqnarray} 
 &&
\omega_1=\frac{( ( {-\Delta})^{n/4}P^T) P^T+P^T ( {-\Delta})^{n/4}P^T-
 ( {-\Delta})^{n/4}(P^TP^T)}{2}\,\label{omega1}\\[5mm]
&&\omega_2={( ( {-\Delta})^{n/4}P^T )P^N+P^T ( {-\Delta})^{n/4}P^N-
 ( {-\Delta})^{n/4}(P^TP^N)},\label{omega2}\\[5mm]
&&\omega=\frac{( ( {-\Delta})^{n/4} P^T) P^T- P^T ( {-\Delta})^{n/4} P^T}{2}\,.\label{Omega}
\end{eqnarray}

We observe that Theorem \ref{comm2}    implies     that 
  $\omega_1,\omega_2 \in L^{(2,1)}(\R)\, $. Moreover it holds
  $$
  \|\omega_1\|_{L^{(2,1)}},\|\omega_2\|_{L^{(2,1)}} \lesssim\|P^T\|^2_{\dot H^{n/2}}\,.$$
The matrix  $\omega$ is   {\bf antisymmetric}.\par
\medskip
{\bf Proof of Proposition  \ref{EulEq}.}\par
From Propositions \ref{pr-I.1} and \ref{pr-I.2} it follows that $u$ satisfies in a weak sense the equations  \rec{euler1} and \rec{eqstructbis}.\par
The key point is to estimate the 
 the
  terms $( ( {-\Delta})^{n/4} P^T) ( {-\Delta})^{n/4} u$ and rewrite equation \rec{eqstructbis} in a different way.\par\medskip

{\bf $\bullet$  Re-writing of $( ( {-\Delta})^{n/4} P^T) ( {-\Delta})^{n/4} u$\,.}\par\medskip
  
\begin{eqnarray*}
( ( {-\Delta})^{n/4} P^T) ( {-\Delta})^{n/4} u&=&( ( {-\Delta})^{n/4} P^T)(P^T ( {-\Delta})^{n/4} u+P^N ( {-\Delta})^{n/4} u)\\[5mm]
&=&
[( ( {-\Delta})^{n/4} P^T) P^T][P^T ( {-\Delta})^{n/4} u)]\\[5mm]
&& +[( ( {-\Delta})^{n/4} P^T) P^N][P^N ( {-\Delta})^{n/4} u)]\,.
\end{eqnarray*}
Now we have
\begin{equation}\label{pt}
( ( {-\Delta})^{n/4} P^T) P^T=\omega_1+\omega+\frac{ ( {-\Delta})^{n/4}P^T}{2}\,;
\end{equation}

 and 

 \begin{eqnarray}\label{pn}
 ( ( {-\Delta})^{n/4} P^T) P^N&=&(  ( {-\Delta})^{n/4} P^T) P^N+P^T ( {-\Delta})^{n/4} P^N\nonumber\\[5mm]
 && - ( {-\Delta})^{n/4}(P^TP^N)-P^T ( {-\Delta})^{n/4}P^N\nonumber\\  &=&\omega_2+P^T ( {-\Delta})^{n/4} P^T\\[5mm]
 &=& \omega_2+\omega_1-\omega+\frac{ ( {-\Delta})^{n/4}P^T}{2}\,.\nonumber
\end{eqnarray}

 Thus
 
 \begin{eqnarray}
 \frac{( ( {-\Delta})^{n/4}P^T) (P^T \Delta ^{n/4} u)}{2}&=&\omega_1 ( P^T \Delta ^{n/4} u) +\omega ( P^T \Delta ^{n/4} u) \label{TT} \\ [5mm]
 \frac{( ( {-\Delta})^{n/4} P^T) (P^N \Delta ^{n/4} u)}{2}&=&(\omega_1+\omega_2) (P^N \Delta ^{n/4} u) -\omega (P^N \Delta ^{n/4} u)\,. \label{TN}  \end{eqnarray}

{\bf $\bullet$   Re-writing of equation \rec{eqstructbis}}.
Equation \rec{eqstructbis} can be re-written as follows:
\begin{eqnarray}\label{eqstruct3}
&& (  ( {-\Delta})^{n/4} (P^N  ( {-\Delta})^{n/4}u)=  (( {-\Delta})^{n/4} \bar {\cal{R}}  )f(P^N,u)\nonumber\\[5mm]
  &&+
 \underbrace{  ( {-\Delta})^{n/4} [P^T  ( {-\Delta})^{n/4}u]- ( {-\Delta})^{\frac{n-1}{2}}   [P^T  ( {-\Delta})^{1/2}u]}_{(3)}\nonumber\\[5mm]
 &&
 +  ( {-\Delta})^{\frac{n-1}{2}}  [P^T  ( {-\Delta})^{1/2}u]- ( {-\Delta})^{n/4} [P^T ( {-\Delta})^{n/4}u]+P^T( {-\Delta})^{\frac{n}{2}}  u\nonumber\\[5mm]
 &&\quad\pm ( ( {-\Delta})^{n/4} P^T) ( ( {-\Delta})^{n/4} u)\\[5mm]
&&= \underbrace{ (  ( {-\Delta})^{\frac{n-1}{2}}  )(P^T  ( {-\Delta})^{1/2}u)- T_n(P^T,u)}_{(4)}\underbrace{-( ( {-\Delta})^{n/4} P^N) ( ( {-\Delta})^{n/4} u)}_{(5)}\,.\nonumber
 \end{eqnarray}
 
  The  term $(3)$ in \rec{eqstruct3} is in
$W^{-n/2, (2,\infty) }\,$ by   Corollary \ref{interp} . The term $(4)$ is in $\dot W^{-n/2,(2,\infty)}(\R^n)$ by Theorem \ref{comm1}  and    Corollary \ref{interp}\,. 
We finally observe that in $(4)$ we use the fact that $P^T( {-\Delta})^{\frac{n}{2}}  u=0$ and in $(5)$ the fact that $ ( {-\Delta})^{n/4} P^T=- ( {-\Delta})^{n/4} P^N\,.$ \par

\par
\medskip

Given $u,Q$ we set
\begin{eqnarray*} 
R(Q,u)&=&(  ( {-\Delta})^{n/4}  )(Q  ( {-\Delta})^{n/4}u)-(  ( {-\Delta})^{\frac{n-1}{2}}  )(Q  ( {-\Delta})^{1/2}u)\\[5mm]
&+&
 (  ( {-\Delta})^{\frac{n-1}{2}}  )(Q  ( {-\Delta})^{1/2}u)-  T_n(Q,u)
 \,.
\end{eqnarray*}
 
We remark that $R(P^T,u)$ is the sum of $(3), (4) $ in \rec{eqstruct3}\,.\par\medskip
 {\bf $\bullet$ Re-writing of $( ( {-\Delta})^{n/4} P^N)  ( {-\Delta})^{n/4} u$\,.}
\par\medskip
 We have
\begin{eqnarray*}
( ( {-\Delta})^{n/4} P^N) ( {-\Delta})^{n/4} u&=& (( {-\Delta})^{n/4}P^N)(P^T( ( {-\Delta})^{n/4} u)+P^N( ( {-\Delta})^{n/4} u)))\,.
\end{eqnarray*}
 \par\medskip
 We estimate
$( ( {-\Delta})^{n/4} P^N) P^T( ( {-\Delta})^{n/4} u)$ and $(  ( {-\Delta})^{n/4} P^N) P^N( ( {-\Delta})^{n/4}u)$\,.
We have
\begin{eqnarray*}
( ( {-\Delta})^{n/4} P^N) P^T&=&-(  ( {-\Delta})^{n/4} P^T) P^T\\[5mm]
&=&-\omega_1-\omega -\frac{( ( {-\Delta})^{n/4} P^T)}{2}\\[5mm]
&=&
-\omega_1-\omega +\frac{( ( {-\Delta})^{n/4} P^N)}{2}\,,
\end{eqnarray*}
and 
\begin{eqnarray*}
 (  ( {-\Delta})^{n/4} P^N) P^N&=&-(  ( {-\Delta})^{n/4} P^T) P^N\pm   P^T(( {-\Delta})^{n/4}P^N )\\[5mm]
&=&
-[(  ( {-\Delta})^{n/4} P^T) P^N+P^T(  ( {-\Delta})^{n/4} P^N) -  ( {-\Delta})^{n/4}(P^N P^T) ]\\[5mm]
&+&
P^T( ( {-\Delta})^{n/4} P^N)\\[5mm]
&=&-\omega_2-P^T( ( {-\Delta})^{n/4} P^T)\\[5mm]
&=&
\mbox{by \rec{TT}}\\[5mm]
&=&-\omega_2-\omega_1+\omega +\frac{( ( {-\Delta})^{n/4} P^N)}{2}\,.
\end{eqnarray*}
Thus
\begin{eqnarray}
\frac{(  ( {-\Delta})^{n/4} P^N) P^T ( {-\Delta})^{n/4} u}{2}&=&-\omega_1(P^T ( {-\Delta})^{n/4} u)-\omega (P^T ( {-\Delta})^{n/4} u)\label{NT} \\[5mm]
\frac{(  ( {-\Delta})^{n/4} P^N) P^N ( {-\Delta})^{n/4} u}{2}&=&-\omega_2(P^N ( {-\Delta})^{n/4} u)-\omega_1(P^N ( {-\Delta})^{n/4} u)\label{NN}\\ 
& &+\omega (P^N ( {-\Delta})^{n/4} u) \nonumber \,.
\end{eqnarray}

By combining \rec{TT}, \rec{TN}, \rec{NT} , \rec{NN}  we obtain

\begin{eqnarray}\label{Eulerbis}
 ( {-\Delta})^{n/4}\left(\begin{array}{l}
P^T ( {-\Delta})^{n/4} u\\ 
 P^N ( {-\Delta})^{n/4} u\end{array}\right)&=&  2\tilde \Omega_1\left(\begin{array}{l}
P^T ( {-\Delta})^{n/4} u\\ 
 P^N ( {-\Delta})^{n/4} u\end{array}\right)+\tilde \Omega_2\\ 
&+&
2\left(\begin{array}{cc}
- \omega &  \omega \\ 
 \omega  & -  \omega  \end{array}\right)
\left(\begin{array}{l}
P^T ( {-\Delta})^{n/4} u\\
 P^N ( {-\Delta})^{n/4} u\end{array}\right)\,,\nonumber
\end{eqnarray}
 where
 $\tilde \Omega_1 $ and  $\tilde  \Omega_2 $ are given by
 $$
\tilde \Omega_1=\left(\begin{array}{cc}
-\omega_1 & -(\omega_1+\omega_2) \\
\omega_1 &  (\omega_1+\omega_2)\end{array}\right)\,;$$

$$
\tilde \Omega_2=\left(\begin{array}{c}
 T_n(P^T,u)\\
 R(P^T,u)+\bar{\cal{R}} ( {-\Delta})^{n/4}f(P^N,u)
\end{array}\right)\,.$$
  
The matrix
$$
\Omega=2\left(\begin{array}{cc}
- \omega &  \omega \\
 \omega  & -   \omega \end{array}\right)$$
is antisymmetric\,.
 
We observe that from the estimates on the operators $T_n$, $R$ and $f$   it follows that
   $\tilde \Omega_2\in \dot {W}^{-n/2, (2,\infty) }(\R^n,\R^{2m})$ and
 \begin{equation}\label{tildeomega1}
\|\tilde\Omega_2\|_{\dot {W}^{-n/2, (2,\infty) }(\R^n,\R^{2m})}\lesssim\left( \|P^N\|_{\dot H^{n/2}(\R^n)}+ \|P^T\|_{\dot H^{n/2}(\R^n)}\right)\,\| ( {-\Delta})^{n/4} u\|_{L^{(2,\infty)}(\R^n)}\,.
\end{equation}
  On the other hand  $\tilde  \Omega_1\in L^{(2,1)}(\R^n,{\cal{M}}_{2m\times 2m})$ and
   \begin{equation}\label{tildeomega2}
\|\tilde\Omega_1\|_{ L^{(2,1)}(\R^n,{\cal{M}}_{2m\times 2m})}\lesssim ( \|P^N\|^2_{\dot H^{n/2}(\R^n)}+ \|P^T\|^2_{\dot H^{n/2}(\R^n)})\,.
\end{equation}

~~\hfill$\Box$

\par\bigskip

Now we prove  {Theorem} \ref{regn2hm}.
\par\medskip
{\bf Proof of Theorem \ref{regn2hm}.} 
From Proposition  \ref{EulEq}  it follows that $$v=(P^T( ( {-\Delta})^{n/4} u), P^N ( ({-\Delta})^{n/4} u)$$ solves equation 
    \rec{Eulerbis}  which  of the type \rec{eqschbis} up to the term $\tilde\Omega_1$ and $\tilde\Omega_2$.

  We aim at obtaining that  $ ( {-\Delta})^{n/4} u\in L^{p}_{loc}(\R)$, for all $p\ge 1\,.$    To this purpose  we take $\rho>0$ such that
$$\|\Omega\|_{L^2(B(0,\rho)}, \|P^T\|_{\dot H^{n/2}(B(0,\rho))}, \|P^N\|_{\dot H^{n/2}(B(0,\rho))}\le \varepsilon_0,$$ with $\varepsilon_0>0$ small enough. 
Let $x_0\in B(0,\rho/4)$ and $r\in(0,\rho/8)$. We argue by duality and multiply  \rec{Eulerbis}  by $\phi$ which is given as follows. 
Let $g\in L^{(2,1)}(\R^n)$, with $\|g\|_{L^{(2,1)}}\le 1$ and set $g_{r\alpha}=\11_{B(x_0,r\alpha)}g$, with $0<\alpha<1/4$ and  $\phi= ( {-\Delta})^{-n/4}(g_{r\alpha})\in L^{\infty}(\R^n)\cap \dot W^{n/2,(2,1)}(\R^{n})$\,.
We multiply  both sides of equation \rec{Eulerbis}  by  $\phi$ and we integrate.\par
   By using the same  ``localization arguments" in the proof of Theorem 1.1 and Theorem 1.7 in \cite{DLR2} one can show that
 $v$ satisfies 
  for all
$x_0\in B(0,\rho/4)$ and $0<r<\rho/8\,,$  
 $
 \|v\|_{L^{(2,\infty)}(B(x_0,r))}\le C r^{\beta}\,,$  for some  $\beta\in(0,1/2)\,.$ 
 
By arguing as in Theorem 1.1 in \cite{DLR2}  we deduce that $ v\in L^{p}_{loc}(\R)$, for all $p\ge 1\,.$  Therefore $ ( {-\Delta})^{n/4} u\in L^{p}_{loc}(\R)$, for all $p\ge 1$ as well.

This implies that $u\in C^{0,\alpha}_{loc}$ for some $0<\alpha<1$, since $W^{n/2,p}_{loc}(\R^n)\hookrightarrow C^{0,\alpha}_{loc}(\R^n)$ if $p>2$
(see for instance \cite{AD2}). This concludes the proof.~\hfill$\Box$
 
 \appendix

  \section{ Commutator Estimates}\label{commutators}
  In this appendix we present a series of commutator estimates which have been used in the previous sections. 
 We consider the  Littlewood-Paley   decomposition of unity introduced in Section \ref{defnot}.    
  For every $j\in \Z$ and $f\in{\cal{S}}^{\prime}(\R^n)$  we define the Littlewood-Paley projection operators $P_j$ and $P_{\le j}$ by 
   \begin{eqnarray*}
 \widehat{ P_jf}=\psi_j \hat{f}~~~\widehat{ P_{\le j}f}=\phi_j \hat{f}\,.
 \end{eqnarray*}
 Informally $P_j$ is a frequency projection to the annulus $\{2^{j-1}\le |\xi|\le 2^{j+1}\}$, while
 $P_{\le j}$ is a frequency projection to the ball $\{|\xi|\le 2^{j+1}\}\,.$ We will set
 $f_j=P_j f$ and $f^j=P_{\le j} f$\,.
 
 We observe that  $f^j=\sum_{k=-\infty}^{j} f_k$  and $f=\sum_{k=-\infty}^{+\infty}f_k$ (where the convergence is in ${\cal{S}}^\prime(\R^n)$)\,.\par
 Given $f,g\in {\cal{S}}^\prime(\R^n)$  we can  split  the product in the following way
 \begin{equation}\label{decompbis}
 fg=\Pi_1(f,g)+\Pi_2(f,g)+\Pi_3 (f,g),\end{equation}
 where
 \begin{eqnarray*}
 \Pi_1(f,g)&=& \sum_{-\infty}^{+\infty} f_j\sum_{k\le j-4} g_k= \sum_{-\infty}^{+\infty} f_j g^{j-4}\,;\\
 \Pi_2(f,g)&=& \sum_{-\infty}^{+\infty} f_j\sum_{k\ge j+4} g_k =\sum_{-\infty}^{+\infty} g_j f^{j-4}\,;\\
 \Pi_3(f,g)&=& \sum_{-\infty}^{+\infty}  f_j\sum_{|k-j|< 4} g_k\,.\
 \end{eqnarray*}
 We observe that   for every $j$ we have 
 $$\mbox{supp${\cal{F}}[f^{j-4}g_j]\subset \{2^{j-2}\le |\xi|\le 2^{j+2}\}$};$$
 $$\mbox{supp${\cal{F}}[\sum_{k=j-3}^{j+3}f_jg_k]\subset \{|\xi|\le 2^{j+5}\}$}\,.$$
 The three pieces of the decomposition \rec{decompbis} are examples of paraproducts. Informally the
 first paraproduct $\Pi_1$ is an operator which allows high frequences of $f$ $(\thicksim 2^{nj})$ multiplied by  low frequences of $g$ $(\ll 2^{nj})$ to produce high frequences in the output. The second paraproduct
 $\Pi_2$ multiplies low fequences of $f$ with high frequences of $g$ to produce high fequences in the output. The third paraproduct $\Pi_3$ multiply high frequences of $f$ with high frequences of $g$ to produce comparable or lower frequences in the output. For a presentation of these paraproducts we refer to the reader for instance  to the book \cite{Gra2}\,.
  We recall the following three results whose proof can be found in \cite{DLR2}.
 
 \begin{Lemma}
 For every $f\in {\cal{S}}^{\prime}$ we have
 $$ \sup_{j\in Z}|f^j|\le M(f)\,.$$
 \end{Lemma}
  \begin{Lemma}\label{lemmaprel1}
Let $\psi$ be a Schwartz  radial function such that $supp (\psi)\subset B(0,4)$. Then 
$$
\|\nabla^k {\cal{F}}^{-1}[\psi]\|_{L^1}\le C_{\psi,n} 4^k\,,
$$
where $C_{\psi,n}$ is a positive constant depending on the $C^2$ norm of $\psi$ and the dimension\,.~$\Box$
\end{Lemma}
 \par
\medskip
\begin{Lemma}\label{lemmaprel2}
Let $f\in B^0_{\infty,\infty}(\R^n)$. Then for all $k\in\N$  and for all $j\in Z$ we have 
$$
2^{-kj}\|\nabla^k f_j\|_{L^\infty}\le 4^k \|f_j\|_{L^\infty}\,.~\hfill \Box
$$
\end{Lemma}
 \par
 \bigskip
 In the sequel we suppose that $n>1$ is an odd integer.\par
  
Given $Q\in \dot{H}^{n/2}(\R^n, {\cal{M}}_{\ell\times m}(\R^n))$ $\ell\ge 0$\footnote{ ${\cal{M}}_{\ell\times m}(\R)$ denotes, as usual, the space of $\ell\times m$ real matrices.} and $u\in  \dot H^{n/2}(\R^n,\R^m)$, we introduce the following operators
 \begin{eqnarray}\label{opM1}
M_1(Q,u)&=&\sum_{{{\displaystyle{\mathop{\scriptstyle{1\le |\alpha|\le [n/2]}}_{|\alpha| odd}}}}}\frac{c_{\alpha}}{\alpha !} ( {-\Delta})^{n/4}([ ( {-\Delta})^{n/4-|\alpha|} (\nabla^{\alpha} Q)]\,\nabla^{\alpha} u)\\[5mm]
&+& \sum_{{{\displaystyle{\mathop{\scriptstyle{1\le |\alpha|\le [n/2]}}_{|\alpha| even}}}}}\frac{c_{\alpha}}{\alpha !} ( {-\Delta})^{n/4}([ ( {-\Delta})^{n/4-|\alpha|/2} Q]\ ,\nabla^{\alpha} u)\,;\nonumber
\end{eqnarray}
\begin{eqnarray}\label{opM2}
M_2(Q,u)&=&\sum_{{{\displaystyle{\mathop{\scriptstyle{1\le |\alpha|\le [n/2]}}_{|\alpha| odd}}}}}\frac{c_{\alpha}}{\alpha !} ( {-\Delta})^{n/4}(\nabla^{\alpha} Q\,[ ( {-\Delta})^{n/4-|\alpha|} (\nabla^{\alpha} u)]\ )\\
&+& \sum_{{{\displaystyle{\mathop{\scriptstyle{1\le |\alpha|\le [n/2]}}_{|\alpha| even}}}}}\frac{c_{\alpha}}{\alpha !} ( {-\Delta})^{n/4}(\nabla^{\alpha} Q\,\,[ ( {-\Delta})^{n/4-|\alpha|/2}  u])\nonumber
\end{eqnarray}
where $c_{\alpha}=\partial^{|\alpha|}|x|^{n/2}_{|_{|x|=1}}\,.$ 
\par
\begin{Proposition}\label{M1M2}
 Let $u,Q\in \dot H^{n/2}(\R^n)$. Then $M_1(Q,u),M_2(Q,u)\in \dot W^{-n/2,(2,1)}(\R^n)$ and
 \begin{eqnarray}
  \|M_1(Q,u)\|_{\dot  W^{-n/2, (2,1)}(\R^n)}&\lesssim &||Q||_{\dot{H}^{n/2}(\R^n)}\| ( {-\Delta})^{n/4}u\|_{L^{2}(\R^n)}\,; \label{estM1}\\
    \|M_2(Q,u)\|_{ \dot W^{-n/2,(2,1)}(\R^n)}&\lesssim &||Q||_{\dot{H}^{n/2}(\R^n)}\| ( {-\Delta})^{n/4}u\|_{L^{2}(\R^n)} \,.\label{estM2}
  \end{eqnarray}
  \end{Proposition}
  {\bf Proof of Proposition \ref{M1M2}\,.} We prove \rec{estM1}, being the estimate of \rec{estM2} similar.\par

We recall that for $0<s<n/2$ we have
$$
\dot H^{n/2}(\R^n)\hookrightarrow \dot W^{s,(\frac{n}{s},2)}(\R^n)\,,$$
(see for instance \cite{T1}).\par
Thus  if $n=2p+1>1$, ($p\ge 1$),  is an odd integer number and $0<|\alpha|\le [n/2]$ then 
$$\nabla^{\alpha} u\in L^{(\frac{n}{|\alpha|},2)},~~( ( {-\Delta})^{n/4-|\alpha|/2} (\nabla^{\alpha} Q)),\,( ( {-\Delta})^{n/4-|\alpha/2|} Q)  \in L^{(\frac{n}{n/2-|\alpha|},2)}\,.$$ \par
 
Thus 
by H\"older Inequality the following  products  
$$[ ( {-\Delta})^{n/4-|\alpha|} (\nabla^{\alpha} Q)]\,\nabla^{\alpha} u\,,~~ [ ( {-\Delta})^{n/4-|\alpha/2|} Q]\, \nabla^{\alpha} u ,$$
are in $  L^{(2,1)}(\R^n)$ and
\begin{eqnarray*}
 \|[ ( {-\Delta})^{n/4-|\alpha|} (\nabla^{\alpha} Q)]\,\nabla^{\alpha} u\|_{L^{(2,1)}}&\lesssim & \|[ ( {-\Delta})^{n/4-|\alpha|} (\nabla^{\alpha} Q)]\|_{L^{(\frac{n}{n/2-|\alpha|},2)}}\|\nabla^{\alpha} u\|_{L^{\frac{n}{|\alpha|},2}}
\\
&
\lesssim & \|Q\|_{\dot{H}^{n/2}(\R^n)}\|u\|_{\dot{H}^{n/2}(\R^n)}\,,
\end{eqnarray*}

\begin{eqnarray*}
\|[ ( {-\Delta})^{n/4-|\alpha|/2} Q]\,\nabla^{\alpha} u\|_{L^{(2,1)}}&\lesssim &\|[ ( {-\Delta})^{n/4-|\alpha|/2}   Q]\|_{L^{(\frac{n}{n/2-|\alpha|},2)}}
\,\|\nabla^{\alpha} u\|_{L^{(\frac{n}{|\alpha|},2)}}
\\
&
\lesssim & \|Q\|_{\dot{H}^{n/2}(\R^n)}\|u\|_{\dot{H}^{n/2}(\R^n)}\,.
\end{eqnarray*}

It follows that $M_1(Q,u)\in  \dot W^{-n/2,(2,1)}(\R^n)  $ and \rec{estM1} holds.\,. This concludes the proof of Proposition \ref{M1M2}.~$\Box$
\par
\bigskip
Next we prove  Theorem \ref{comm2}.\par
 \medskip
 {\bf Proof of Theorem \ref{comm2}.}
 We group as follows\,:
 \begin{eqnarray*}
  && \Pi_1[ T_n^*(Q,u)]=
 \underbrace{\Pi_1[ ( {-\Delta})^{n/4}(( ( {-\Delta})^{n/4} Q)\,u)- ( {-\Delta})^{n/2}(Q\, u)]}+\underbrace{\Pi_1[ ( {-\Delta})^{n/4} (  Q\,  ( ( {-\Delta})^{n/4} u))]}\\
 &&
 \Pi_2[ T_n^*(Q,u)]=
 \underbrace{\Pi_2[ ( {-\Delta})^{n/4}(( ( {-\Delta})^{n/4} Q)\,u)]}+\underbrace{\Pi_2[- ( {-\Delta})^{n/2}(Q\, u)+ ( {-\Delta})^{n/4} ( Q\,  ( ( {-\Delta})^{n/4} u)]}\\
&& \Pi_3[ T_n^*(Q,u)]=
 \underbrace{\Pi_3[ ( {-\Delta})^{n/4}(( ( {-\Delta})^{n/4} Q)\,u)]}-\underbrace{\Pi_3[ ( {-\Delta})^{n/2}(Q\, u)]}+\underbrace{\Pi_3[ ( {-\Delta})^{n/4} (  Q\,  ( ( {-\Delta})^{n/4} u))]}\,.
 \end{eqnarray*}
 Some   terms appearing in $ T_n^*(Q,u)$ satisfy a better estimate in the sense that they belong in ${\cal{H}}^1$ or in $\dot{B}^0_{1,1} $. We recall that
  $\dot{B}^0_{1,1}\hookrightarrow {\cal{H}}^1 \hookrightarrow  \dot W^{-n/2,(2,1)}\,.$\par
  
 {\bf  $\bullet$}  Estimate of $\|\Pi_1[ ( {-\Delta})^{n/4}(Q ( {-\Delta})^{n/4}u)]\|_{{\cal{H}}^1(\R^n)}$.
 \begin{eqnarray}\label{P1tn}
 &&
 \|\Pi_1[ ( {-\Delta})^{n/4}(Q ( {-\Delta})^{n/4}u)]\|_{{\cal{H}}^1(\R^n)}\simeq \int_{\R^n}\left(\sum_j \left(2^{\frac{n}{2} j} Q_j ( {-\Delta})^{n/4} u^j)\right)^2\right)^{1/2}dx\\[5mm]
 &&
\lesssim \int_{\R^n}\sup_{j} [( {-\Delta})^{n/4} u^j] \left(\sum_j 2^{nj} Q^2_j\right)^{1/2}dx\nonumber\\[5mm]
 &&
 \lesssim\left(\int_{\R^n}(\sup_{j} ( {-\Delta})^{n/4} u^j)^2 dx\right)^{1/2}\left(\int_{\R^n}\sum_j 2^{nj} Q^2_j\,dx\right)^{1/2}\nonumber\\[5mm]
 &&
 \lesssim
 \|Q\|_{\dot{H}^{n/2}(\R^n)}\|u\|_{\dot{H}^{n/2}(\R^n)}\,.\nonumber
 \end{eqnarray}
 {\bf  $\bullet$}  Estimate of $\|\Pi_3[ ( {-\Delta})^{n/4}(Q ( {-\Delta})^{n/4}u)]\|_{{B_{1,1}^0}(\R^n)}\,.$
 \begin{eqnarray}\label{P3tn}
 &&
 \|\Pi_3[ ( {-\Delta})^{n/4}(Q ( {-\Delta})^{n/4}u)]\|_{\dot B_{1,1}^0(\R^n)}\\[5mm]
&&~~~~ \simeq \sup_{\|h\|_{B^0_{\infty,\infty}\le 1}}\int_{\R^n} \sum_jQ_j ( {-\Delta})^{n/4} u_j)[ ( {-\Delta})^{n/4}h^{j-6}+\sum_{t=j-5}^{j+6} ( {-\Delta})^{n/4} h_t]\nonumber\\[5mm]
 &&
~~~~ \lesssim \sup_{\|h\|_{B^0_{\infty,\infty}\le 1}}
 \|h\|_{B^0_{\infty,\infty}}\int_{\R^n} 2^{n/2 j}|Q_j ||( {-\Delta})^{n/4} u_j|\,dx\nonumber\\[5mm]
 &&~~~~~
 \lesssim
 \|Q\|_{\dot{H}^{n/2}(\R^n)}\|u\|_{\dot{H}^{n/2}(\R^n)}\,.\nonumber
 \end{eqnarray}
 The estimate of $   \Pi_3[ ( {-\Delta})^{n/2}(Q\, u)]$, $\Pi_3[( {-\Delta})^{n/4} (  Q\,  ( ( {-\Delta})^{n/4} u)]$, $\Pi_2[ ( {-\Delta})^{n/4}(( ( {-\Delta})^{n/4} Q)\,u)]$ are similar to \rec{P1tn} and \rec{P3tn} and we omit them.

{\bf  $\bullet$}  Estimate of $\|\Pi_1[  ( {-\Delta})^{n/4}(( ( {-\Delta})^{n/4} Q)\,   u)- ( {-\Delta})^{n/2}(Q\, u)]\|_{\dot  W^{-n/2,(2,1)}(\R^n)} $\,.
  
 \begin{eqnarray}\label{p1comm2}
&&||\Pi_1( ( {-\Delta})^{n/4}( ( {-\Delta})^{n/4} Q) u- ( {-\Delta})^{n/2}(Qu))||_{ W^{-n/2,(2,1)}(\R^n)}\\[5mm]
&&=\sup _{||h||_{ \dot W^{n/2,(2,\infty)}(\R^n)}\le 1} \int_{\R^n}\sum_j\sum _{|t-j|\le 3}[ ( {-\Delta})^{n/4}(  ( {-\Delta})^{n/4} Q_ju^{j-4})- ( {-\Delta})^{n/2}( Q_ju^{j-4})]h_tdx\nonumber\\[5mm]
&& =
\sup _{||h||_{ \dot W^{n/2,(2,\infty)}(\R^n)}\le 1} \int_{\R^n}\sum_j\sum _{|t-j|\le 3} {\cal{F}}[u^{j-4}]
{\cal{F}}[\Delta ^{n/4}Q_j  ( {-\Delta})^{n/4} h_t-Q_j ( {-\Delta})^{n/2}h_t]d\xi\nonumber\\ 
 && =\sup _{||h||_{ \dot W^{n/2,(2,\infty)}(\R^n)}\le 1}\int_{\R^n}\sum_j\sum _{|t-j|\le 3} {\cal{F}}[u^{j-4}](\xi)\nonumber\\[5mm]
&&~~~ \left(\int_{\R^n} {\cal{F}}[Q_j] (\zeta){\cal{F}}[ ( {-\Delta})^{n/4} h_t](\xi-\zeta)( |\zeta|^{n/2}-|\xi-\zeta|^{n/2} ) d\zeta\right) d\xi\,.\nonumber
\end{eqnarray}

 Now we observe that in \rec{p1comm2} we have $|\xi|\le 2^{j-3}$ and $2^{j-2}\le |\eta|\le 2^{j+2}$.
Thus $|\displaystyle\frac{\xi}{\zeta}|\le \frac{1}{2}\,.$ 
Hence
\begin{eqnarray}\label{taylorexpbis}
&& |\zeta|^{n/2}-|\xi-\zeta|^{n/2}=|\xi|^{n/2}\left[1-\left|\frac{\xi}{|\zeta|}-\frac{\zeta}{|\zeta|}\right|^{n/2}\right]\nonumber\\
&&=|\zeta|^{n/2}\left[\sum_{{{\displaystyle{\mathop{\scriptstyle{|\alpha|\ge 1}}_{|\alpha| odd}}}}}\frac{c_\alpha}{\alpha!}\left(\frac{\xi}{|\zeta|}\right)^{\alpha}\left(\frac{\zeta}{|\zeta|}\right)^{\alpha}
+\sum_{{{\displaystyle{\mathop{\scriptstyle{|\alpha|\ge 2}}_{|\alpha| even}}}}}\frac{c_\alpha}{\alpha!}\left(\frac{\xi}{|\zeta|}\right)^{\alpha}\right] \,.
 \nonumber
 \end{eqnarray}
 We may suppose the series in \rec{taylorexpbis} is   convergent if $|\displaystyle\frac{\xi}{|\zeta|}|\le \frac{1}{2}\,,$ otherwise one may consider a
 different Littlewood-Paley decomposition by replacing the exponent $j-4$ with $j-s$, $s>4$ large enough.  \par
 Unlike   the case $n=1$ (see the proof of estimate $(35)$ in \cite{DLR1})   we need to separate two cases:
  $|\alpha|\ge [n/2]+1$ and $1\le |\alpha|\le  [n/2]$\,.
     
  \par
  {\bf Case 1: $ |\alpha|\ge [n/2]+1$\,.}
      Here we use the fact that
 $  \dot W^{n/2,(2,\infty)}(\R^n)\hookrightarrow \dot{B}^0_{\infty,\infty}(\R^n)$ and the crucial property that for every vector field $X\in \dot H^{n/2}(\R^n)$ we have
\begin{eqnarray}\label{equiv}
&&\int_{\R^n}\sum_{j=-\infty}^{+\infty} 2^{-jn}(X^{j})^2 dx
=\int_{\R^n} \sum_{k,\ell}X_kX_{\ell}\sum_{j-4\ge k,j-4\ge \ell}2^{-jn}dx\nonumber\\ [5mm]
&&\simeq \int_{\R^n} \sum_{k}X_k\left(\sum_{|k-\ell|\le 2} X_{\ell} \right)\,2^{-(k-2)n}dx\nonumber\\[5mm]
&&
\mbox{by Cauchy-Schwarz Inequality}\nonumber\\[5mm]
&& \lesssim\int_{\R^n}\left(\sum_{k} 2^{-kn} X_k^2\right)^{1/2}\left(\sum_{k} 2^{-kn} X_k^2\right)^{1/2} dx\\[5mm]
&&=\int_{\R^n}\sum_{j=-\infty}^{+\infty} 2^{-kn}(X_k)^2dx\,,\nonumber\end{eqnarray}
(see also Section 4.4.2 in  \cite{RS}, page 165).\par\medskip
   We are going to estimate
     \begin{eqnarray}\label{p1comm24bis}
  && \sup _{||h||_{\dot W^{n/2,(2,\infty)}}\le 1} 
[  \sum_{{{\displaystyle{\mathop{\scriptstyle{|\alpha|\ge [n/2]+1}}_{|\alpha| odd}}}}}\frac{c_\alpha}{\alpha!}\int_{\R^n}\sum_j  \nabla^{\alpha} u^{j-4}   ( {-\Delta})^{n/4-|\alpha|}(\nabla^{\alpha}Q_j) ( {-\Delta})^{n/4} h_j  dx\nonumber\\[5mm]
  &&~~~~+  \sum_{{{\displaystyle{\mathop{\scriptstyle{|\alpha|\ge [n/2]+1}}_{|\alpha| even}}}}}\frac{c_\alpha}{\alpha!}\int_{\R^n}\sum_j  |\nabla^{\alpha} u^{j-4}  ( {-\Delta})^{n/4-|\alpha|/2}( Q_j) ( {-\Delta})^{n/4} h_j  dx]\,.
  \end{eqnarray}\par
By applying Lemma \ref{lemmaprel2} ($\|( {-\Delta})^{n/4} h_j \|_{\dot{B}^0_{\infty,\infty}}\lesssim 2^{\frac{nj}{2}} 4^{n/2} ||h||_{\dot{B}^0_{\infty,\infty}} $) we get

  \begin{eqnarray}\label{p1comm24}
  && \rec{p1comm24bis} \lesssim \sup _{||h||_{\dot W^{n/2,(2,\infty)}}\le 1} ||h||_{\dot{B}^0_{\infty,\infty}}\nonumber\\
  &&
[  \sum_{{{\displaystyle{\mathop{\scriptstyle{|\alpha|\ge [n/2]+1}}_{|\alpha| odd}}}}}\frac{c_\alpha}{\alpha!}\int_{\R^n}\sum_j 2^{\frac{nj}{2}}|\nabla^{\alpha} u^{j-4}|| ( {-\Delta})^{n/4-|\alpha|}(\nabla^{\alpha}Q_j)| dx\nonumber\\
  &&+  \sum_{{{\displaystyle{\mathop{\scriptstyle{|\alpha|\ge [n/2]+1}}_{|\alpha| even}}}}}\frac{c_\alpha}{\alpha!}\int_{\R^n}\sum_j 2^{\frac{nj}{2}}|\nabla^{\alpha} u^{j-4}|| ( {-\Delta})^{n/4-|\alpha|/2}( Q_j)| dx]\nonumber\\
    &&\lesssim \sup _{||h||_{\dot W^{n/2,(2,\infty)}}\le 1} ||h||_{\dot{B}^0_{\infty,\infty}}\nonumber\\
 && [  \sum_{{{\displaystyle{\mathop{\scriptstyle{|\alpha|\ge [n/2]+1}}_{|\alpha| odd}}}}}\frac{c_\alpha}{\alpha!}2^{2n-4|\alpha|} \\
&&  \left(\int_{\R^n}\sum_j 2^{(n-2|\alpha|)(j-4)}|\nabla^{\alpha} u^{j-4}|^2dx\right)^{1/2} \left(\int_{\R^n}\sum_j 2^{2|\alpha|j}| ( {-\Delta})^{n/4-|\alpha|}(\nabla^{\alpha}Q_j)|^2dx\right)^{1/2}]\nonumber\\
 &&
+
  [ \sum_{{{\displaystyle{\mathop{\scriptstyle{|\alpha|\ge [n/2]+1}}_{|\alpha| even}}}}}\frac{c_\alpha}{\alpha!}2^{2n-4|\alpha|}\nonumber\\
&&  \left(\int_{\R^n}\sum_j 2^{(n-2|\alpha|)(j-4)}|\nabla^{\alpha} u^{j-4}|^2\,dx\right)^{1/2} \left(\int_{\R^n}\sum_j 2^{2|\alpha|j}| ( {-\Delta})^{n/4-|\alpha|/2}( Q_j)|^2dx\right)^{1/2}] \nonumber\\
&&
  \lesssim \|Q\|_{\dot H^{n/2}(\R^n)}\|u\|_{\dot H^{n/2}(\R^n)}\nonumber\,.
  \end{eqnarray}
  \par
  \medskip
  {\bf Case 2: $ 1\le |\alpha|\le [n/2]$\,.}
  In this case we apply Proposition \ref{M1M2}. \par
 We have:
  \begin{eqnarray*}
 && \sup _{||h||_{\dot W^{n/2,(2,\infty)}}\le 1}  \\
  &&
[   \sum_{{{\displaystyle{\mathop{\scriptstyle{1\le |\alpha|\le [n/2]}}_{|\alpha| odd}}}}}\frac{c_\alpha}{\alpha!}\int_{\R^n}\sum_j  \nabla^{\alpha} u^{j-4}   ( {-\Delta})^{n/4-|\alpha|}(\nabla^{\alpha}Q_j) ( {-\Delta})^{n/4} h_j  dx\nonumber\\[5mm]
  &&+   \sum_{{{\displaystyle{\mathop{\scriptstyle{1\le |\alpha|\le [n/2]}}_{|\alpha| even}}}}}\frac{c_\alpha}{\alpha!}
  \int_{\R^n}\sum_j  |\nabla^{\alpha} u^{j-4}  ( {-\Delta})^{n/4-|\alpha|/2}( Q_j) ( {-\Delta})^{n/4} h_j  dx]\nonumber\\[5mm]
  &&
  \lesssim \sup _{||h||_{\dot W^{n/2,(2,\infty)}}\le 1} \| ( {-\Delta})^{n/4} h\|_{L^{(2,\infty)}}\nonumber\\[5mm]&&
[  \sum_{{{\displaystyle{\mathop{\scriptstyle{1\le |\alpha|\le [n/2]}}_{|\alpha| odd}}}}}\frac{c_\alpha}{\alpha!}
    \| ( {-\Delta})^{n/4-|\alpha|}(\nabla^{\alpha}Q)\|_{\dot W^{|\alpha|,(\frac{n}{n/2-|\alpha|},2)}}
    \|\nabla^{\alpha} u\|_{\dot W^{\frac{n}{2}-|\alpha|,(\frac{n}{|\alpha|},2)}}\nonumber\\[5mm]
    &&
    +
    [ \sum_{{{\displaystyle{\mathop{\scriptstyle{1\le |\alpha|\le [n/2]}}_{|\alpha| even}}}}}\frac{c_\alpha}{\alpha!}
    \| ( {-\Delta})^{n/4-|\alpha|/2} Q\|_{\dot W^{|\alpha|,(\frac{n}{n/2-|\alpha|},2)}}
    \|\nabla^{\alpha} u\|_{\dot W^{\frac{n}{2}-|\alpha|,(\frac{n}{|\alpha|},2)}}\nonumber\\[5mm]
    &&
  \lesssim
  \sup _{||h||_{\dot W^{n/2,(2,\infty)}}\le 1}
  \| ( {-\Delta})^{n/4} h\|_{L^{(2,\infty)}} \|Q\|_{\dot H^{n/2}(\R^n)}\|u\|_{\dot H^{n/2}(\R^n)}\nonumber\\[5mm]
  &&
  \lesssim \|Q\|_{\dot H^{n/2}(\R^n)}\|u\|_{\dot H^{n/2}(\R^n)}\,.\nonumber
  \end{eqnarray*}

The estimate of $\Pi_2[  ( {-\Delta})^{n/4}( Q\,  ( ( {-\Delta})^{n/4}  u))- ( {-\Delta})^{n/2}(Q\, u) ]$ is analogous to \rec{p1comm2}  and we omit it. This concludes the proof of Theorem \ref{comm2}\,.~\hfill$\Box$
\par
\bigskip
  \par
\vskip0.5truecm
The next result permits us to estimate the right hand side of equation \rec{eqstructbis}\,.\par
   We denote by $r^{\prime}$ the coniugate of $1<r<+\infty$. \par
\begin{Theorem}\label{l2bis}
 Let $n>2$, $1<r<\frac{2n}{n-2},$   $h\in L^{r^{\prime}}(\R^n)$, $Q\in\dot H^{n/2}(\R^n)$. Then
\begin{equation}\label{l2est2}
  ( {-\Delta})^{ \frac{n}{4}-\frac{1}{2}}(Qh)- Q ( {-\Delta})^{ \frac{n}{4}-\frac{1}{2}} h \in \dot W^{-(\frac{n}{2}-1),r^{\prime}}(\R^n)\,,
\end{equation}
and
\begin{equation}\label{l2est4}
\| ( {-\Delta})^{ \frac{n}{4}-\frac{1}{2}}(Qh)- Q ( {-\Delta})^{ \frac{n}{4}-\frac{1}{2}} h\|_{\dot W^{-(\frac{n}{2}-1),r^{\prime}}(\R^n)}\lesssim \|h\|_{L^r}\|Q\|_{\dot H^{n/2}(\R^n)}\,.
\end{equation}
\end{Theorem}
  Theorem \ref{l2bis} implies ``by duality"  the following result.
\begin{Theorem}\label{l2}
 Let $n>2$, $1<r<\frac{2n}{n-2},$ 
     $Q\in\dot H^{n/2}(\R^n)$, $f\in\dot W^{ \frac{n}{2}-1,r}(\R^n)$ then
\begin{equation}\label{l2est}
Q ( {-\Delta})^{ \frac{n}{4}-\frac{1}{2}} f- ( {-\Delta})^{ \frac{n}{4}-\frac{1}{2}}(Qf) \in L^r(\R^n)\,,
\end{equation}
and
\begin{equation}\label{l2est1}
\|Q ( {-\Delta})^{ \frac{n}{4}-\frac{1}{2}} f- ( {-\Delta})^{ \frac{n}{4}-\frac{1}{2}}(Qf)\|_{L^r(\R^n)}\lesssim \|Q\|_{\dot H^{n/2}(\R^n)}\|f\|_{\dot W^{ \frac{n}{2}-1,r}}\,.~~~\Box
\end{equation}
\end{Theorem}

\par
\vskip0.5truecm

{\bf Proof of Theorem \ref{l2bis}\,.}
Throughout the proof we use  the following embeddings:
$$
\dot W^{\frac{n}{2}-1,r}(\R^n)\hookrightarrow L^s(\R^n),\quad \frac{1}{s}=\frac{1}{r}-\frac{\frac{n}{2}-1}{n}\,;
$$
$$
\dot H^{n/2}(\R^n)\hookrightarrow \dot W^{\frac{n}{2}-1,(q,2)}(\R^n)\quad \frac{1}{q}=\frac{1}{2}-\frac{1}{n}=\frac{\frac{n}{2}-1}{n}\,.
$$
We also use the fact that
\begin{equation}\label{magic}
\frac{1}{r^{\prime}}+\frac{1}{s}+\frac{1}{q}=1\,.
\end{equation}
{\bf $\bullet$ Estimate of $\|\Pi_1[ ( {-\Delta})^{ \frac{n}{4}-\frac{1}{2}}(Qh)]\|_{\dot W^{-(\frac{n}{2}-1),r^{\prime}}(\R^n)}\,.$}\par
\begin{eqnarray}\label{p11}
&&
\|\Pi_1[ ( {-\Delta})^{ \frac{n}{4}-\frac{1}{2}}(Qh)]\|_{\dot W^{-(\frac{n}{2}-1),r^{\prime}}(\R^n)}\nonumber\\[5mm]
&&~~\simeq
\sup_{\|g\|_{\dot W^{(\frac{n}{2}-1),r }(\R^n)}\le 1}\int_{\R^n}\sum_j Q_j h^j  ( {-\Delta})^{ \frac{n}{4}-\frac{1}{2}} g_j dx\nonumber\\[5mm]
&&
~~\lesssim
\sup_{\|g\|_{\dot W^{(\frac{n}{2}-1),r }(\R^n)}\le 1}
\int_{\R^n}\sup_{j} h^j\sum_j 2^{(\frac{n}{2}-1)j} Q_j 2^{-(\frac{n}{2}-1)j} ( {-\Delta})^{\frac{n}{4}-\frac{1}{2}} g_j\,dx\\[5mm]
&&
~~\mbox{\bf by generalized H\"older Inequality}\nonumber\\[5mm]
&&\lesssim \sup_{\|g\|_{\dot W^{(\frac{n}{2}-1),r }(\R^n)}\le 1}\|h\|_{L^{r^{\prime}}}\| ( {-\Delta})^{ \frac{n}{4}-\frac{1}{2}}Q\|_{L^q}\|g\|_{L^s}\nonumber\\[5mm]
&&
~~\lesssim 
\|h\|_{L^{r^{\prime}}}\|Q\|_{\dot H^{n/2}(\R^n)}\,.\nonumber
\end{eqnarray}

{\bf $\bullet$ Estimate of $\|\Pi_3[ ( {-\Delta})^{ \frac{n}{4}-\frac{1}{2}}(Qh)]\|_{\dot W^{-(\frac{n}{2}-1),r^{\prime}}(\R^n)}\,.$}\par
\begin{eqnarray}\label{p31}
&&
\|\Pi_3[ ( {-\Delta})^{ \frac{n}{4}-\frac{1}{2}}(Qh)]\|_{\dot W^{-(\frac{n}{2}-1),r^{\prime}}(\R^n)}\nonumber\\[5mm]
&&\simeq
\sup_{\|g\|_{\dot W^{(\frac{n}{2}-1),r }(\R^n)}\le 1}\int_{\R^n}\sum_j Q_j h_j  ( {-\Delta})^{ \frac{n}{4}-\frac{1}{2}} g^j dx\nonumber\\[5mm]
&&
\lesssim
\sup_{\|g\|_{\dot W^{(\frac{n}{2}-1),r }(\R^n)}\le 1}
\int_{\R^n}\sup_{j} 2^{-(\frac{n}{2}-1)j}( ( {-\Delta})^{ \frac{n}{4}-\frac{1}{2}} g^j)\sum_j 2^{(\frac{n}{2}-1)j} Q_j  h_j\,dx\\[5mm]
&&
\mbox{\bf by generalized H\"older Inequality}\nonumber\\[5mm]
&&\lesssim \sup_{\|g\|_{\dot W^{(\frac{n}{2}-1),r }(\R^n)}\le 1}\|h\|_{L^{r^{\prime}}}\| ( {-\Delta})^{ \frac{n}{4}-\frac{1}{2}}Q\|_{L^q}\|g\|_{L^s}\nonumber\\[5mm]
&&
\lesssim 
\|h\|_{L^{r^{\prime}}}\|Q\|_{\dot H^{n/2}(\R^n)}\,.\nonumber
\end{eqnarray}
\par
The estimates of $\Pi_1[Q ( {-\Delta})^{ \frac{n}{4}-\frac{1}{2}}h]$ and $\Pi_3[Q ( {-\Delta})^{ \frac{n}{4}-\frac{1}{2}}h]$ are similar to \rec{p11} and \rec{p31} and we omit them.\par

{\bf $\bullet$ Estimate of $\|\Pi_2[ ( {-\Delta})^{ \frac{n}{4}-\frac{1}{2}}(Qh)-Q ( {-\Delta})^{ \frac{n}{4}-\frac{1}{2}}h]\|_{\dot W^{-(\frac{n}{2}-1),r^{\prime}}(\R^n)}\,.$}\par
We denote by $\tilde c_{\alpha}$ the coefficients of the Taylor expansion of 
$|x|^{\frac{n}{2}-1}$ at $x=1\,.$\par
\begin{eqnarray}\label{provv3}
&&\|\Pi_2[ ( {-\Delta})^{ \frac{n}{4}-\frac{1}{2}}(Qh)-Q ( {-\Delta})^{ \frac{n}{4}-\frac{1}{2}}h]\|_{\dot W^{-(\frac{n}{2}-1),r^{\prime}}(\R^n)}\nonumber\\[5mm]
&&
\simeq 
\sup_{\|g\|_{\dot W^{(\frac{n}{2}-1),r }(\R^n)}\le 1}
\int_{\R^n}\left( ( {-\Delta})^{ \frac{n}{4}-\frac{1}{2}}(Q^jh_j)-Q^j ( {-\Delta})^{ \frac{n}{4}-\frac{1}{2}}h_j\right)g_j\,dx\,.
\end{eqnarray} 
Now we argue as 
  in \rec{p1comm2} and we get\par
\begin{eqnarray}\label{provv4}
&&\rec{provv3}
\lesssim \sup_{\|g\|_{\dot W^{(\frac{n}{2}-1),r }(\R^n)}\le 1}
[   \sum_{{{\displaystyle{\mathop{\scriptstyle{1\le |\alpha| }}_{|\alpha| odd}}}}}\frac{\tilde c_\alpha}{\alpha!}\int_{\R^n}\sum_j  \nabla^{\alpha} Q^{j-4}   ( {-\Delta})^{ \frac{n}{4}-\frac{1}{2}-|\alpha|}(\nabla^{\alpha}h_j)  g_j  dx\nonumber\\[5mm]
  &&+   \sum_{{{\displaystyle{\mathop{\scriptstyle{1\le |\alpha| }}_{|\alpha| even}}}}}\frac{\tilde c_\alpha}{\alpha!}
  \int_{\R^n}\sum_j  |\nabla^{\alpha} Q^{j-4}  ( {-\Delta})^{ \frac{n}{4}-\frac{1}{2}-|\alpha|/2}( h_j)  g_j  dx]\nonumber\\[5mm]
&&
\mbox{by Lemma \ref{lemmaprel2}}\nonumber\\[5mm]
&&
\lesssim \sup_{\|g\|_{\dot W^{(\frac{n}{2}-1),r }(\R^n)}\le 1}\|Q\|_{B^0_{\infty,\infty}}\nonumber\\[5mm]
&&
[\sum_{{{\displaystyle{\mathop{\scriptstyle{1\le |\alpha| }}}}}}\frac{\tilde c_\alpha}{\alpha!}2^{-2|\alpha|}\int_{\R^n} 2^{|\alpha|j}2^{-(\frac{n}{2}-1)j} |( {-\Delta})^{ \frac{n}{4}-\frac{1}{2}-|\alpha|/2}(h_j)|
2^{(\frac{n}{2}-1)j}|g_j |dx]\nonumber\\[5mm]
&&
\lesssim
  \sup_{\|g\|_{\dot W^{(\frac{n}{2}-1),r }(\R^n)}\le 1}\|Q\|_{B^0_{\infty,\infty}}\|h\|_{L^{r^{\prime}}}\| ( {-\Delta})^{ \frac{n}{4}-\frac{1}{2}}g\|_{L^r}\nonumber\\[5mm]
  &&
  \lesssim
  \|Q\|_{\dot H^{n/2}}\|h\|_{L^{r^{\prime}}}\,.\nonumber\,~~\hfill\Box
  \end{eqnarray}
\par
\medskip
 Since $\frac{2n}{n-2}>2$ we can  now apply  the interpolation Theorem 3.3.3 in \cite{Hel} and obtain the following:
 \begin{Corollary}\label{interp}
 Let $n>2$ ,
     $Q\in\dot H^{n/2}(\R^n)$, $f\in\dot H^{ \frac{n}{2}-1}(\R^n)$ then
\begin{equation}\label{l2est2bis}
Q ( {-\Delta})^{ \frac{n}{4}-\frac{1}{2}} f- ( {-\Delta})^{ \frac{n}{4}-\frac{1}{2}}(Qf) \in L^{(2,\infty) }(\R^n)\,,
\end{equation}
and
\begin{equation}\label{l2est3}
\|Q ( {-\Delta})^{ \frac{n}{4}-\frac{1}{2}} f- ( {-\Delta})^{ \frac{n}{4}-\frac{1}{2}}(Qf)\|_{L^{(2,\infty) }(\R^n)}\lesssim \|Q\|_{\dot H^{n/2}(\R^n)}\|f\|_{\dot W^{ \frac{n}{2}-1,(2,\infty)}}\,.
\end{equation}
\end{Corollary}
 
   \medskip
   
\vskip0.5cm
We finally recall a commutator estimate obtained in \cite{CRW} .
\begin{Lemma}\label{crw}
Let $p\ge 1$,  $Q\in BMO(\R^n)$, $u\in L^p(\R^n)$ and let ${\cal {P}}$\footnote{We recall that a pseudo-differential operator ${\cal{P}}$ can be formally defined as
 $$
 {\cal{F}}[{\cal{P}}f(x)]= \sigma (x,\xi){\cal{F}}[f],
 $$ where $\sigma$, the symbol of ${\cal{P}}$, is a complex-valued function defined $\R^n\times\R^n$. If $\sigma(x,\xi)=m(\xi)$ is independent of $x$, then ${\cal{P}}$ is the Fourier  multiplier associated with $m$\,. Given $k\in\Z$ we say that $\sigma$ is of order $k$ if 
   for every multi-indexes $\beta,\alpha\in \N^n$
 $$ |D^{\beta}_xD^\alpha_{\xi} \sigma(x,\xi)|\le C_{\alpha,\beta} |\xi|^{k-\alpha}\,.$$} a pseudo-differential operator of order zero. Then
${\cal{P}}(Qu)-Q{\cal{P}}u\in L^p(\R^n)$
and
$$
\|{\cal{P}}(Qu)-Q{\cal{P}}u\|_{L^p(\R^n)}\lesssim\|Q\|_{BMO(\R^n)}\|u\|_{L^p(\R^n)}\,.~~\Box
$$
\end{Lemma}
The interpolation Theorem 3.3.3 in \cite{Hel}, and  Lemma  \ref{crw} imply the following result.
\begin{Corollary}\label{crwbis}
Let $Q\in BMO(\R^n)$, $u\in L^{(2,\infty) }(\R^n)$ and let ${\cal {P}}$ a pseudo-differential operator of order zero. Then
${\cal{P}}(Qu)-Q{\cal{P}}u\in L^{2,\infty }(\R^n)$
and
$$
\|{\cal{P}}(Qu)-Q{\cal{P}}u\|_{ L^{(2,\infty) }(\R^n)}\lesssim\|Q\|_{BMO(\R^n)}\|u\|_{ L^{(2,\infty)}(\R^n)}\,~~\Box
$$
\end{Corollary}
We    observe that Corollary \ref{crwbis} implies  that for every $h\in L^{(2,1)}(\R^n), u\in L^{(2,\infty) }(\R^n)$
the operator $u{{\cal{P}}} h-({\cal{P}} u) h\in {\cal{H}}^1(\R^n)$ and
\begin{equation}\label{estPs}
\|u{{\cal{P}}} h-({\cal{P}} u) h\|_{{\cal{H}}^1(\R^n)}\lesssim \|u\|_{ L^{(2,\infty)}(\R^n)}\|h\|_{ L^{(2,1)}(\R^n)}\,.
\end{equation}\par

\par
\bigskip
{\bf Aknowledgment.} The author would like to thank Tristan Rivi\`ere for helpful comments.

  \end{document}